\documentclass{article}
\usepackage{amsthm}  
\usepackage{amssymb}
\usepackage{amsmath}
\usepackage{mathptmx}

\newcommand{\mbb}[1]{\mathbf{#1}}

\newcommand{\mrm}[1]{\mathrm{#1}}

\newcommand{\Q}{\mbb{Q}}

\newcommand{\fn}[3]{{#1:#2\rightarrow #3}}

\theoremstyle{definition}
\newtheorem{defn}{Definition}[section]
\newtheorem{rmk}[defn]{Remark}

\theoremstyle{plain}
\newtheorem{thm}[defn]{Theorem}

\newtheorem{prop}[defn]{Proposition}

\begin{document}

\title{On the properness of the moduli space of stable surfaces}
\author{Michael A. van Opstall}
\date{}

\maketitle

\begin{abstract}
The moduli spaces of stable surfaces serve as compactifications of the moduli
spaces of canonical models of smooth surfaces in the same way the moduli 
spaces of stable curves compactify the moduli spaces of smooth curves. However,
the natural definition of the moduli functor of stable surfaces allows in
extra components parameterizing surfaces which do not smooth. This article
verifies that these components are also proper.
\end{abstract}

\section{Introduction}

The theory of moduli spaces of stable curves has clearly shown the usefulness
of a compactification of a moduli space of algebraic varieties which is itself
a moduli space. By allowing curves with the mildest possible singularities 
into the moduli problem, one obtains a tractable projective variety which
compactifies the space of all smooth curves of genus $g$. 

A similar space exists which provides a projective compactification of the 
components of the moduli space of canonically polarized surfaces (projective
surfaces with ample canonical class and at most rational double points). The
varieties parameterized by the boundary of this space have mild (in some sense)
singularities, but singularities which are more complex than mere normal
crossings points. The compactification is called the moduli space of stable
surfaces.

Projectivity of this space follows from properness due to a result of
Koll\'ar (\cite{k:pcm}, Theorem 4.12). Properness follows from the 
boundedness results of
Alexeev (\cite{al:bk2}, Section 7) and verification of the valuative criterion of properness,
which in turn follows from the existence of semi-stable canonical models
(cf. Chapter 7 of \cite{km:bgav}). 

After allowing limits of canonically polarized surfaces into the moduli 
problem, their deformations must also be admitted. In the case of curves,
nodal curves all have unobstructed deformations, so this adds no new components
to the moduli space. On the other hand, there are degenerations
of smooth surfaces having deformations to singular surfaces which admit
no smoothings. If these extra components of the moduli space are excluded
somehow from the moduli functor, it is not clear how to define a scheme 
structure on the moduli space. What is the tangent space at these bad points
if we artificially throw away components? Note that in his thesis 
\cite{hack:th}, Hacking gives a functorial way ({\em smoothable families})
of removing these components in the case of
log surfaces occurring as degenerations of the plane marked with a smooth
curve.

One solution to this problem is not to throw away these components. In this
case, a problem lingers: are these new components projective? The 
aforementioned results of Koll\'ar and Alexeev reduce this question to the
valuative criterion of properness. Checking this criterion is the subject
of this article. This question was suggested by J. Koll\'ar at the AIM
workshop ``Compact moduli spaces and birational geometry''. I am very
grateful to Brendan Hassett for pointing out some errors in an older
version.

\subsection{Sketch of the argument}

The basic method of finding ``canonical'' limits of families of smooth
varieties is to apply the semistable reduction theorem (a good statement is
Theorem 7.17 of \cite{km:bgav}) to obtain (possibly
after a finite base change) a family of reduced simple normal crossing 
varieties completing the original family. The ambiguity of choices of the
special fibers is then removed by taking the relative canonical model of
this family. For curves, this amounts to blowing down rational curves in
the fibers with self-intersection -1 or -2 (that is, rational curves in
a fiber meeting the rest of the fiber in one or two points). For surfaces, 
the process is significantly more complex, and uses essentially all the
operations of the minimal model program.

When the general member of the family is singular, the total space may be
non-normal, so the machinery of semistable reduction and canonical models
does not directly apply. One could try to extend this machinery to the case
of varieties with at worst normal crossings singularities in codimension one,
but there are difficulties. It is better to normalize the original family,
work on the components one by one, and attempt to reassemble the desired 
result.

For curves this is straightforward: pull the irreducible components of the
total space apart, marking them with their (horizontal) conductors. Take
a semistable resolution of the component pieces which also desingularizes
(hence separates components of) the conductors. Now take the relative
log canonical model of the pair of the total space of the family marked with
its conductors. The conductors of the new family are nonsingular curves
birational to the original conductors, so everything glues back together
in the end. The separated conductors do not meet in the relative canonical
model, since a relative canonical model is an {\em lc morphism}. Of course, 
since all nodal curves are smoothable, there is no
need to check the valuative criterion on families of curves whose members
are all singular: the nonsingular curves are an open dense set in the moduli
space.

There are a few problems in extending this argument to higher dimensions. 
First,
taking relative canonical models does not recover the fibers of the original
family in general. For example, a surface which is smooth except for a simple
elliptic singularity and which has ample canonical class is a stable surface
whose minimal desingularization is also stable. Each of these surfaces is a
log canonical model. Starting with a family of such singular surfaces and
applying the algorithm above results in a family whose general member is
the minimal desingularization of the general member of the original family.

The second problem is that in higher dimensions, it is not obvious that the
conductors will glue back together, especially after making the additional
modifications necessary to overcome the first problem. However, the fact that
they do can be seen as a consequence of the separatedness of the moduli space
of stable varieties one dimension lower. Also, in general, components of 
the conductor may meet after taking the relative canonical model, but the
conductors may still be glued together to admissible varieties (in the surface
case, this leads to {\em degenerate cusps}).

\subsection{Notation}

I will use the notations and basic definitions from higher dimensional
geometry following \cite{km:bgav} and \cite{k:fa}. 
To keep notation simpler, I will use 
some conventions. First, if $(X,D)$ is a pair obtained from $(Y,B)$ by some 
birational morphism from $X$ to $Y$ or from $Y$ to $X$, $D$ will be the
birational (``strict'' or ``proper'') transform of $B$ on $X$ (obtained by pushing forward $B$ by the
birational morphism in question or its inverse). I will write a morphism
$\pi:(X,D)\rightarrow C$ with a pair as domain to emphasize that canonical
models to be taken are log canonical models. Consequently, I will frequently
drop the adjectives ``relative'' and ``log'' and speak simply of canonical
models. If $C^0$ is an open set of $C$, then $X^0$ is the part of $X$ lying
over this open set, and similarly for $D^0$. All restrictions of a morphism
$\pi$ to smaller sets will remain denoted $\pi$. Finally, all equivalences
are equivalences of $\Q$-divisors.

\section{Preliminaries}

The most useful definition of semi-log canonical for the purposes of this
article is:

\begin{defn}
A variety $X$ has {\em semi-log canonical (slc)} singularities if
\begin{enumerate}
\item $X$ is $S_2$;
\item the singularities of $X$ in codimension one are (double) normal 
crossings;
\item $X$ is $\Q$-Gorenstein, i.e. $\omega_X^{[N]}$ (the reflexive hull of
the $N$th tensor power of $\omega_X$) is locally free for some $N$;
\item the pair $(X^\nu, D)$ consisting of the normalization of $X$ marked with
its conductor is log canonical (lc).
\end{enumerate}
A {\em stable variety} is a projective variety with slc singularities such
that $\omega_X^{[N]}$ is an ample invertible sheaf for large and divisible $N$.
\end{defn}

We will use this definition in both directions: normalizing an slc variety
produces a collection of lc pairs, and gluing a collection of lc pairs along
their boundaries yields an slc variety.

Families of stable varieties are not represented by a separated space, so an
additional condition is necessary.

\begin{defn}
Let $\pi:X\rightarrow B$ be a flat projective morphism whose fibers are 
stable varieties.
\begin{enumerate}
\item $\pi$ is {\em weakly $\Q$-Gorenstein} if $X$ is $\Q$-Gorenstein.
\item $\pi$ is {\em $\Q$-Gorenstein} if $X$ is $\Q$-Gorenstein and the
reflexive powers $\omega_{X/B}^{[n]}$ of the relative dualizing sheaf 
commute with arbitrary base change.
\end{enumerate}
\end{defn}

The second definition is better, since it leads to a natural deformation
theory, but a fundamental remaining 
question in the theory is whether these notions differ. For the purposes
of this article they do not: 

\begin{prop}
Notation as in the definition. If the base $B$ is a smooth curve and the 
fibers are curves or surfaces, then weakly
$\Q$-Gorenstein implies strongly $\Q$-Gorenstein. 
\end{prop}

\begin{proof}
This is proved in \cite{hack:th}, Proposition 10.14 if the general fiber
is canonical. The main point in the proof is to show a certain divisor is 
$S_2$. If the general fiber is canonical, then the total space is canonical, 
which allows Hacking to conclude that the divisor in question is actually 
Cohen-Macaulay. In general, families with smooth base and slc fibers have
slc total spaces (assuming Inversion of Adjunction), which allows us to 
conclude in the same manner as Hacking that his divisor $Z$ is $S_2$.
\end{proof}

Therefore, I 
will prove the result for the weaker condition, and drop the adverb ``weakly''.

The following is Definition 7.1 in {\em loc. cit.}:

\begin{defn}
A nonconstant morphism $\fn{f}{(X,D)}{C}$ to a smooth curve with $X$ normal
and $D$ an effective $\Q$-divisor is called {\em lc} if for every closed point
$c\in C$, $(X,D+f^{-1}(c))$ is an lc pair.
\end{defn}

\begin{prop}\label{fiboflc}
If $\fn{\pi}{(X,D)}{C}$ is an lc morphism and $c\in C$ is a closed
point, 
then 
$(f^{-1}(c),D|_{f^{-1}(c)})$ is slc if $D$ is $S_2$.
\end{prop}

\begin{proof}
Cf. Lemma 7.4(2) of \cite{km:bgav}. The boundaries of lc pairs are normal
crossings in codimension one.
\end{proof}

\begin{rmk}
The final clause of the proposition is a real one. If $D$ is the boundary
divisor in a family of stable varieties, there is no guarantee that $D$
should be $S_2$. For families of log surfaces this is true by the results
of Hassett's article \cite{hass:cm}. 
\end{rmk}

Due to the last remark, I restrict attention to families of surfaces below.
The results will be more generally valid if Hassett's result mentioned in
the remark can be extended.

\section{Normal generic fibers}

First I will show that $\Q$-Gorenstein families of stable lc pairs can be 
completed to families of stable pairs. For the remainder of this section,
fix the following notation: $\pi:(X,D)\rightarrow C$ is a flat, projective 
morphism to (the germ of) a smooth curve, and there is a nonempty open subset 
$C^0\subset C$ such that:
\begin{enumerate}
\item $K_{X^0}+D^0$ is $\Q$-Cartier and $\pi$-ample;
\item the fibers of $\pi$ over $C^0$ are lc pairs;
\end{enumerate}

\begin{thm}\label{lcextend}
There exists a finite and surjective 
base change $C'\rightarrow C$ such that the pullback of $(X^0,D^0)\rightarrow
C^0$ extends to an lc morphism $(X',D')\rightarrow C'$ with $K_{X'}+D'$
relatively ample. 
\end{thm}

\begin{proof}
A key word in the statement of this theorem is {\em extends}.

Suppressing the base change in notation, assume $(X,D)\rightarrow C$ admits
semistable resolution, and without loss of generality, $C\backslash C^0$ is a 
single (closed) point $0\in C$.

Suppose $\fn{f_i}{(X_i,D_i)}{(X,D)}$ for $i=1,2$ are semistable resolutions
of $(X,D)$, and write:
\[
K_{X_i^0}+D_i^0=f_i^*(K_X^0+D^0)+\sum a^{(i)}_j E^{(i)}_j
\]
where without loss of generality all $E$ dominate $C^0$. Denote the closure
of the $E$ in $X_i$ by the same symbol.
Denote by $(X_i^c, D_i^c)$ the relative canonical models (over $(X,D)$) of
$(X_i,D_i-\sum_{b^{(i)}_k<0}a^{(i)}_kE^{(i)}_k)$ and the corresponding
morphisms to $(X,D)$ by $g_i$. I make two claims:
\begin{enumerate}
\item $(X_i^c,D_i^c)$ is independent of $i$;
\item the fibers of $\pi\circ g_i$ over $C^0$ agree with the fibers of
$\pi$.
\end{enumerate}

\begin{proof}[Proof of Claim 1]
Let $(\tilde{X},\tilde{D})$ be a semistable resolution dominating
$X_1$ and $X_2$ by morphisms $h_1$ and $h_2$. Write
\[
K_{\tilde{X}}+\tilde{D}=h_i^*(K_{X_i}+D_i)+\sum c_i G_i.
\]
Then since $X_i$ is smooth, the $c_i$ are all positive. Therefore applying
Corollary 3.53 from \cite{km:bgav}, we conclude that the relative log 
canonical model
of 
\[
(\tilde{X},\tilde{D}-\mrm{exc.~divisors~with~negative~discrepancy})
\]
agrees with that of
\[
(X_i,D_i-\mrm{exc.~divisors~with~negative~discrepancy})
\]
\end{proof}

Now denote by $(X',D')$ the common relative canonical model $(X_i^c,D_i^c)$.

\begin{proof}[Proof of Claim 2]
Let $c\in C^0$. The morphism $\fn{g}{(X'_c,D'_c)}{(X_c,D_c)}$ is a morphism
onto a canonical model. Write
\[
K_{X'_c}+D_c'=g^*(K_{X_c}+D_c)+\sum a_k E_k
\]
following the form above. 
\[
K_{X'_c}+D_c'-\sum_{a_k<0} a_k E_k-g^*(K_{X_c}+D_c)
\]
is effective and exceptional, so the canonical model of $(X'_c,D'_c)$
coincides with $(X_c,D_c)$ ({\em ibid.}).
\end{proof}
\end{proof}

\section{General families}

In this section, I will first cover the case of families of surfaces, where
the conjectures necessary for the existence of the moduli space are theorems,
and then I will dicuss the general case.

The following theorem proves the properness of the irreducible components of
the moduli of stable surfaces.

\begin{thm}[Main Theorem]\label{main}
Suppose $\pi:X\rightarrow C$ is a flat projective morphism to the germ of a
smooth curve and that $C^0$ is a nonempty open subset of $C$ such that
\begin{enumerate}
\item $X^0$ is $\Q$-Gorenstein;
\item the fibers of $\pi$ over $C^0$ are slc surfaces;
\item $K_{X^0}$ is $\pi$-ample.
\end{enumerate}
Then there exists a finite base change $C'\rightarrow C$ and a $\Q$-Gorenstein
family $\tilde{X}$ of stable surfaces extending the pullback of $X|_{C^0}$ to
$C'$.
\end{thm}

\begin{proof}
Denote by $X^\nu=\coprod (X_i,D_i)$ the normalization of $X$ marked by 
its horizontal conductors. I will ignore any vertical conductors (throwing
them away makes $(X_i,D_i)$ ``even more'' lc) and write
the prime decomposition $D_i=\sum_j D_{ij}$. Take a base change such that
all the pairs $(X_i,D_i)$ admit semistable resolution. Apply Theorem 
\ref{lcextend} to these pairs. If some $D_{ij}\cong D_{kl}$ as a result of
the way $X$ is glued together, making the base change does not affect
this isomorphism. Furthermore, over $C^0$ (base change suppressed in notation),
the birational transforms of $D_{ij}$ and $D_{kl}$ are isomorphic (by Claim 2
in the proof of \ref{lcextend}), and over
all of $C$, these birational transforms are families of nodal curves, since
the boundaries of log canonical surface pairs are nodal curves.

Denote a semistable resolution of the pair $(X_i,D_i)$ by 
$(\tilde{X_i},\tilde{D}_i)$. Denote the relative log canonical model of
$(\tilde{X_i},\tilde{D}_i)$ with respect to $K_{\tilde{X}_i}+\tilde{D}_i$ 
by $(\bar{X}_i,\bar{D}_i)$
and
the relative canonical model with respect to 
\[
K_{\tilde{X}_i}+\tilde{D}_i-\sum_{a_i<0}a_iE_i
\]
(as in the proof of \ref{lcextend}) by $(\hat{X}_i,\hat{D}_i)$. Since the 
fibers of $(\bar{X}_i,\bar{D}_i)$ as well as those of $(\hat{X}_i,\hat{D}_i)$
are lc pairs, $\bar{D}_i$ and $\hat{D}_i$ are families of nodal curves. Since
$\tilde{X}_i$ and $\tilde{D}_i$ are smooth, by adjunction $K_{\bar{D}_i}=
K_{\tilde{X}_i}+\tilde{D}_i|_{\bar{D}_i}$. Therefore, $\bar{D}_i$ is the
relative canonical model of a family of nodal curves (with $K_{\bar{D}_i}$
big), so it is a family of stable curves (cf. Proposition 3.3 of
\cite{hass:lsr}). $\hat{D}_i$ is a family of nodal
curves dominated by the family $\bar{D}_i$, so it too is a family of stable
curves. Therefore, the limiting
curve $\hat{D}_{i,0}$ is uniquely determined. At this point 

The morphisms to each of these log canonical models may only contract some or
all of the $E$, or exceptional divisors mapping to $0\in C$. Restricting to a 
fiber, the $E$ restrict to curves which are
not components of any $D$. The collapsing of $E$ may result in components of
the conductor coming together. If this is so, then two familes of stable
curves are glued together, resulting in a family of stable curves, and the
``mates'' of these two components of the conductor must come together over
the general fiber (again, by Claim 2 of \ref{lcextend}). Two components of
the conductor cannot meet as a result of collapsing divisors in the central
fiber, since then the pair $(\bar{X},\bar{D}+X_0)$ would not be lc, 
so the relative canonical model would not be an lc morphism.
Therefore, the families $\hat{D}_i$ and $\bar{D}_i$
coincide (possibly another base change is needed to ensure this, since two
families of stable curves can have the same fibers and not be the same
family), if we renumber the $D_i$ to refer to {\em connected} components
of the conductor rather than {\em irreducible} components. 

Therefore, the identifications among various conductors are preserved through
the process of taking base change, semistable resolution, and relative 
canonical models. We may glue the various $(\hat{X}_i,\hat{D}_i)$ together
to obtain a new family which is the $\bar{X}$ in the statement of the theorem.

Since all of the components $\hat{X}_i$ and the identified loci $\hat{D}_i$
are uniquely determined by the original family, the limit is unique.
\end{proof}

As noted above, this theorem together with the work of Alexeev and Koll\'ar
proves the projectivity of the connected components of the moduli spaces
of stable surfaces.

In general, what is needed for the existence of the moduli spaces and their
projectivity is:
\begin{enumerate}
\item Local closedness of the moduli functor - this would follow from results
of Hassett and Kovacs \cite{hk:lc} for the weakly $\Q$-Gorenstein functor
and from results of Hacking \cite{hack:th} or unpublished results of
Abramovich and Hassett for the strongly $\Q$-Gorenstein functor, if we knew
that having semi-log canonical singularities is an open condition. This in 
turn follows from the minimal model program: see \cite{karu:mmpbd}, Lemma 2.6.
A recent preprint \cite{kaw:ia} of Kawakita proves the necessary Inversion of
Adjunction-type result without recourse to the minimal model program.
\item Boundedness of the moduli functor - what is needed here is a 
generalization of Alexeev's result. Specifically, one needs to know that the
set of all log canonical pairs $(X,D)$ such that $K_X+D$ is ample and with 
fixed $(K_X+D)^n$ is bounded. For {\em smoothable} stable surfaces, Karu 
({\em loc. cit.}, Theorem 1.1) shows
that the boundedness follows from the MMP (in one dimension higher than
the moduli problem under consideration), but this is not true in general.
\item Existence of relative canonical models of semistable resolutions -
the proofs in this article make it clear why this is necessary.
\item Cohen-Macaulayness of limits of Cohen-Macaulay varieties. See the remark
in the first section.
\end{enumerate}

With all of these results in place, 
the argument at the end of \ref{main} that the ``extra''
exceptional divisors to be blown down meet the conductor in points 
is not valid in higher dimensions, but this
should not matter. The family $\bar{D}$ of
stable varieties dominates the family $\hat{D}$ of slc varieties, hence
$\hat{D}$ is also a family of stable varieties and has a unique limit. 

One should also consider, following Alexeev, moduli spaces of stable pairs
$(X,D)$. It is essential in higher dimensions that $D$ be reduced (or at least
that some condition be imposed). In this case, the arguments given here
go through.

\end{document}